\documentclass[12pt,a4paper,english,psamsfonts,reqno]{amsart}
\usepackage{amsfonts}
\usepackage{amsmath}
\usepackage{amssymb}
\usepackage{bbm}
\usepackage{mathrsfs}
\usepackage[dvips]{graphicx}

\usepackage{amsthm}

\theoremstyle{plain}\newtheorem{lemma}{Lemma}
\newtheorem{proposition}[lemma]{Proposition}
\newtheorem{theorem}[lemma]{Theorem}
\newtheorem{corollary}[lemma]{Corollary}

\theoremstyle{definition}
\newtheorem{example}[lemma]{Example}
\newcommand{\N}{{\bf N}}

\newcommand{\C}{{\bf C}}

\newcommand{\rme}{{\rm e}}
\newcommand{\rmd}{{\rm d}}
\newcommand{\cA}{\mathcal{A}}

\newcommand{\cE}{\mathcal{E}}

\newcommand{\cJ}{\mathcal{J}}

\newcommand{\sig}{\sigma}
\newcommand{\alp}{\alpha}
\newcommand{\bet}{\beta}
\newcommand{\gam}{\gamma}

\newcommand{\lam}{\lambda}
\newcommand{\del}{\delta}
\newcommand{\eps}{\varepsilon}

\newcommand{\Spec}{{\rm Spec}}

\newcommand{\norm}{\Vert}
\renewcommand{\Re}{{\rm Re}\;}
\renewcommand{\Im}{{\rm Im}\;}

\newcommand{\Note}{\underbar{Note}{\hskip 0.1in}}
\newcommand{\hash}{\#}

\newcommand{\pr}{\prime}
\newcommand{\Epm}{E_{-+}}
\usepackage{hyperref}
\title{Perturbations of Jordan Matrices}
\author{E B Davies}
\address{
Department of Mathematics\\
King's College\\
Strand\\
London WC2R 2LS} \email{E.Brian.Davies@kcl.ac.uk}
\author{Mildred Hager}
\address{CMLS\\Ecole polytechnique\\91128 Palaiseau
  C{\'e}dex\\France\\UMR 7640}
\email{hager@math.polytechnique.fr}
\date{16 November 2006}
\begin{document}

\begin{abstract} We consider perturbations of a large Jordan
matrix, either random and small in norm or of small rank.
In both cases we show that most of the eigenvalues of the
perturbed matrix are very close to a circle with centre at
the origin. In the case of random perturbations we obtain
an estimate of the number of eigenvalues that are well
inside the circle in a certain asymptotic regime. In the
case of finite rank perturbations we completely determine
the spectral asymptotics as the size of the matrix
increases. The paper provides an elementary illustration
of some standard techniques of spectral theory.
\end{abstract}

\maketitle

\section{Introduction\label{sect1}}
It is well known that the eigenvalues of large non-normal
matrices can be highly unstable under very small
perturbations. In this note we discuss a very simple
example of this phenomenon. We show that in a wide variety
of cases almost all of the eigenvalues of a slightly
perturbed Jordan matrix lie near a circle with centre at
the origin, with high probability in the random case. We
also examine the exceptional eigenvalues, which remain
well inside the circle.

A quantitative measure of spectral instability is provided
by the notion of pseudospectra, which become interesting
when the operator involved is far from being normal; see
\cite{TE,Dav2} for detailed discussions and many
references. If $\del >0$ the $\delta$-pseudospectra of an
operator $A$ are defined by
\begin{equation} \label{pseudo}
\begin{split} \Spec_{\delta}(A)& = \Spec(A)
\cup \{ z \notin \Spec(A): \|(A-z)^{-1}\| > \delta^{-1} \}
\\ & = \bigcup_{\{K: \|K\|< 1\}} \Spec(A+\delta K) \,,
\end{split}
\end{equation}
where $\Spec$ denotes the spectrum of a matrix. The second
equality in \eqref{pseudo} implies that a perturbation of
$A$ of size $\delta$ can move the eigenvalues anywhere
inside $\Spec_{\delta}(A)$. In particular the computed
eigenvalues of a large matrix may be very inaccurate if
$\Spec_{\delta}(A)$ is a large region, where $\delta$ is
the rounding error of the computations. In this note, we
study this phenomenon in some detail for the Jordan block
matrix, perturbed either by a matrix of small rank, in
which case the analysis is much sharper, or by a random
matrix with a small norm. The problem studied in this
paper was proposed by Zworski, who showed how the general
methods of Sj\"ostrand and Zworski (\cite{sjzwgru}) could
be adapted to this particular setting. Our results go
beyond the theory of Lidskii (\cite{Lid,MBO}) by allowing
larger (but still extremely small) perturbations, for
which the Puiseux series is not convergent and the
eigenvalues are not where the first few terms of that
series would predict.

We define the standard $N\times N$ Jordan matrix $J$ by
\[
J_{r,s}:=\left\{ \begin{array}{ll} 1&\mbox{ if $s=r+1$,}\\
0&\mbox{otherwise,} \end{array}\right.
\]
where $r,s=1,...,N$, and we always assume $N>2$ from now
on. Figure 1 shows the results of a MATLAB computation of
$\Spec(J +\delta K)$, when $N = 100$, $\delta =10^{-10}$
and $K$ is a complex gaussian random matrix. Our goal is
to explain the form of the figure and others obtained by
similar methods.

\begin{figure}
\begin{center}
\includegraphics[width=10cm]{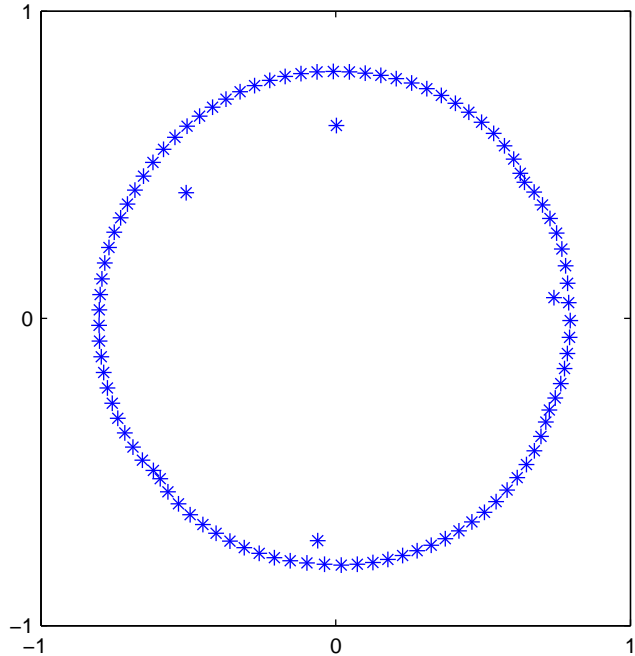}
\end{center}
\caption{Eigenvalues of a small random perturbation of the
Jordan matrix $J$, with $ N = 100 $ and $ \delta =
10^{-10} $.} \label{fig}
\end{figure}

Most of the eigenvalues accumulate around a circle with
centre at the origin, hence far away from $\Spec(J)=\{
0\}$, even though the perturbation is very small in norm
(with a high probability). In Section~\ref{hag}, we
explain the origin of this instability and prove that it
is very likely to happen for this type of perturbation.
The section illustrates the methods and ideas of \cite{Ha}, \cite{hasj}
in a very concrete setting.

If one adds a strictly upper triangular matrix to $J$ then
the spectrum is not changed. In Section \ref{dav}, we
therefore concentrate on perturbations whose non-zero
entries are all close to the bottom left-hand corner of
the matrix; some generalizations are considered in
Section~\ref{general}. We give a complete asymptotic
analysis of the spectrum for all such perturbations. The
problem is described in detail in the following section.
The equation to be solved is written down in
Theorem~\ref{setup}. The asymptotic form of the solutions
to this equation is described in Theorem~\ref{main1}, that
will also be related to a Grushin problem as in
\cite{sjzwgru}, see also Section~\ref{hag}, and then in
much more detail in Theorem~\ref{main2}.
Section~\ref{sectpencil} describes the asymptotics of the
zeros of some more general polynomials.

\section{Random perturbations}
\label{hag} Let us start by investigating the
pseudospectra of the Jordan matrix. Let $D(w,r)= \{z
\in\C:  |z-w|\leq r\}$.
\begin{lemma} \label{pseudospec}
If $0<|z|<1$, then
\begin{equation} \label{resnorm}
|z|^{-N} \leq \|(J-z)^{-1}\| \leq N|z|^{-N},
\end{equation}
which implies that for $\delta <1/N$,
\begin{equation} \label{pseudospecjor}
D(0,\delta^{1/N}) \subseteq \Spec_{\delta}(J)
\subset D(0,(\delta N)^{1/N}) \ .
\end{equation}
\end{lemma}

\begin{proof}
Let $ e(z)=(1,...,z^{N-1})^{\pr}$. If $0<|z|<1$, the
identity
\[
\|(J-z)e(z)\|=|z|^{N}\leq |z|^N\| e(z)\|
\]
implies the lower bound in (\ref{resnorm}). We have the
following expression for the resolvent for $z \neq 0$:
\[
(J-z)^{-1}_{r,s} =
 \left\{ \begin{array}{ll}
                     -z^{r-s-1} & \mbox{if $1\leq r \leq s \leq N$,}\\
            0    & \mbox{otherwise.}
                     \end{array}
              \right.
\]
We conclude that for $|z| <1$,
\begin{equation}
\|(J-z)^{-1}u\|^{2} \leq (\sum_{j=1}^{N} |z|^{-j})^{2}
\|u\|^{2} \leq (N|z|^{-N})^{2} \|u\|^{2} \ .
\end{equation}
\end{proof}

Using (\ref{pseudo}), we deduce that the eigenvalues of a
perturbation of norm $\delta$ may be anywhere in
$D(0,(\delta N)^{1/N})$. In Figure 1 most, but not all, of
the eigenvalues are close to the boundary of this disc,
and our goal is to understand why this is the case.

\begin{theorem} \label{th1}
Let $K$ be a $N \times N$ random matrix such that
\begin{equation} \label{normest}
P[\|K\| < 1] \geq 1 - p_{1}(N) \ ,
\end{equation}
and
\begin{equation} \label{lb}
P[|K_{N,1}|<s ] \leq p_{2}(s;N)  \ .
\end{equation}
Then for any $0<\delta<\frac{1}{2N}$, $\alpha \geq
\delta$, and $\sigma>0$, with probability at least
\begin{equation}
1-p_{1}(N) - p_{2}(5\alpha ;N) \ ,
\end{equation}
we have
\begin{equation} \label{specbound}
\Spec(J+\del K)\subseteq D(0,(\delta N)^{1/N}) \ ,
\end{equation}
and
\begin{equation} \label{specest}
\# (\Spec(J+ \delta K) \cap D(0,(\delta
N)^{1/N}e^{-\sigma}) ) \leq \frac{1}{\sigma}(\ln N
- \ln \alpha) \ .
\end{equation}
\end{theorem}

The following theorem is obtained by putting
$\alpha=N^{-3}$ and estimating $p_{1}(N)$ and
$p_{3}(5\alpha, N)$. We see that for any fixed $\sigma
>0$, the proportion of eigenvalues that lie in the annulus
\[
\{ z: (\delta N)^{1/N}e^{-\sigma} \leq |z| \leq (\delta
N)^{1/N} \}
\]
converges to one with probability one as $N \to \infty$.

\begin{theorem} \label{cor1}
Let $\tilde{K}$ be a $N\times N$ random matrix with its
entries independently and identically distributed
according to a complex gaussian law centered at $0$ and of
variance $1$. Let $K=\tilde{K}/N^{2}$. Then for any $0<
\delta \leq N^{-3}$ and $\sigma>0$, with probability at
least $1-26/N^{2}$, (\ref{specbound}) is valid and
\begin{equation}
\# \big( \Spec(J+ \delta K) \cap D(0,(\delta
N)^{1/N}e^{-\sigma}) \big) \leq
\frac{4}{\sigma}\ln N \ .
\end{equation}
\end{theorem}

Choosing $\delta=\gamma^{N}$, we also obtain the following
result.
\begin{corollary} \label{cor2}
Let $K$ be as in Theorem~\ref{cor1}. Then for any $0<
\gamma \leq N^{-3/N}$ and $\sigma>0$, with probability at
least $1- 26 N^{2}$,
\begin{equation}
\Spec(J+\gamma^{N} K)\subseteq D(0,\gamma N^{1/N})
\end{equation}
and
\begin{equation}
\# \big( \Spec(J+ \gamma^{N} K) \cap D(0,\gamma
N^{1/N} e^{-\sigma}) \big) \leq
\frac{4}{\sigma}\ln N \ .
\end{equation}
\end{corollary}

\subsection{The Grushin problem}
\label{grushin} For the proof of Theorem \ref{th1}, we set
up a Grushin problem as in \cite[Sect. 2.2]{sjzwgru}. Let
$A \in M_{N}(\C)$, let $m \in \N$, and let $R_{+}$ and
$R_{-}'$ be $m\times N$ matrices. We put
\begin{equation} \label{gru}
\cA  = \left( \begin{array}{ccc}
                     A & R_{-} \\
                     R_{+} & 0
                     \end{array}
              \right) \in M_{N+m}(\C) \ .
\end{equation}
\begin{lemma}(Schur, Grushin)\label{Schur}
If $\cA$ is invertible, and
\[
\cE  = \left( \begin{array}{ccc}
                     E & E_+ \\
                     E_- & \Epm
                     \end{array}
              \right)
\]
is the matrix inverse of $\cA$, then $A$ is invertible if
and only if  $\det(\Epm)\not= 0$.
\end{lemma}
\begin{proof}
If $\cA$ is invertible and $\det(\Epm)\not= 0$, then
\begin{equation} \label{inverse}
A(E-E_{+}\Epm^{-1}E_{-}) = 1 \ ,
\end{equation}
hence $A$ is invertible. The converse affirmation goes
along the same path.
\end{proof}
\begin{corollary} \label{r+-}
If
\begin{equation} \label{grujo}
m=1 \ ,  \ R_{-}=e_{N} \ ,  \ R_{+}=e_{1}' \ ,
\end{equation}
where $e_{1},e_{2},...,e_{N}$ is the standard basis of
column vectors in $\C^{N}$, then $\cA$ is invertible if
and only if $\det(\tilde{A})\not= 0$, where $\tilde{A}$ is
obtained by deleting the first column and last row of $A$.
In that case
\[
\Epm=(-1)^N\frac{\det(A)}{\det(\tilde{A})}.
\]
\end{corollary}

\begin{example} \label{grujor} If $A=J-zI$ and $R_{\pm}$ are as in (\ref{grujo}), we write $\cJ=\cA$. By the previous corollary $\cJ$ is invertible, and
\[
\cE_{r,s} =
 \left\{ \begin{array}{ll}
                     z^{r-s-1} & \mbox{if $s+1\leq r\leq N+1$,}\\
                    z^{r-1} & \mbox{if $1\leq r\leq N+1$ and $s=N+1$,} \\
            0    & \mbox{otherwise.}
                     \end{array}
              \right.
\]
Assuming $|z|\leq 1$ we deduce that
\begin{equation}\label{Ebounds} \begin{split}
& \|\cE(z)\| \leq N+1,\hspace{8mm} \|E(z)\| \leq
N,\hspace{8mm}
\|E_\pm(z)\| \leq N^{1/2},\hspace{8mm} \Epm(z)=z^N. \\
& \|\cE(0)\| \leq 1,\hspace{8mm} \|E(0)\| \leq
1,\hspace{8mm} \|E_\pm(0)\| \leq 1,\hspace{8mm} \Epm(0)=0.
\end{split}
\end{equation}
Finally, using (\ref{inverse}), we can also find the
explicit expression for the resolvent.
\end{example}

\subsection{Perturbation}
Let us assume that $\|K\|<1$. Then by
(\ref{pseudospecjor}), (\ref{specbound}) holds. We analyze
the part of the spectrum within $D(0,(\del
N)^{1/N})$ in more detail by using the Grushin
problem for $J+\delta K$. We will show that the matrix
$\mathcal{J}^{\delta}=\mathcal{A}$ obtained by putting
$A=J+\delta K-zI$ in (\ref{gru}) may be inverted by using
a Neumann series. Denoting the inverse by $\cE^\del$, this
implies that $ \Spec(J+\del K)\cap D(0,R)$ coincides with
the set of zeros of $\Epm^\del (z)$ such that $|z|\leq R$.

\begin{lemma}\label{Epm}
Let $\|K\|<1$, and let $\delta <1/2N$. Then if $N\geq 2$,
for any $R<1$,
\begin{equation}
\Spec(J+\del K) \cap D(0,R) = \{ z \in D(0,R) ;
\Epm^{\del}(z)=0 \} \ ,
\end{equation}
where
\begin{equation}
\Epm^\del(z)=z^N - \del p_K(z) + q_{\del K}(z)
\label{eqEpm}
\end{equation}
and
\begin{equation} \label{p}
p_{K}(z) = \sum_{r,s=0}^{N-1} K_{(N-r),s+1} z^{r+s} \ .
\end{equation}
If $|z| < (\delta N)^{1/N}=R$, we have
\begin{equation} \label{estepmgen}
\| \Epm^\del \|_\infty \leq 3\del N \ ,
\end{equation}
in $L^{\infty}(D(0,R))$, and
\begin{equation}
| \Epm^\del (0) | \geq \del( | K_{N,1}| - 2 \delta) \ .
\end{equation}
\end{lemma}

\begin{proof} The resolvent expansion implies invertibility
provided $\norm \del K\|  \| \cE\norm<1$, which is proved
by using (\ref{Ebounds}) and $\del \|K\|N <1/2$. Moreover
we have the following Neumann series expansion for the
inverse:
\begin{align}
\mathcal{E}^{\delta} & = \mathcal{E}^{0} + \left(
\begin{array}{ccc} \sum_{j\geq 1} E(-\delta K E)^{j} &
 \sum_{j\geq 1} (-E \delta K )^{j} E_
{+}) \\
\sum_{j\geq 1} E_{-}(-\delta K E)^{j} & \sum_{j\geq 1}
E_{-}(-\delta K E)^{j-1}(-\delta K E_{+}) \end{array}
\right) \ .
\end{align}
Evaluating the bottom right coefficient of each term
yields (\ref{eqEpm}) with
\[
\begin{split}
p_K(z)&=E_-KE_+ \ , \\
q_{\del K}(z)&= \del^2 E_- K E K E_+ - \del^3 E_- K E K E
K E_+ + ... \ .
\end{split}
\]
The bounds on the various quantities now follow by
combining (\ref{Ebounds}), $\|K\|<1$, $\delta < 1/2N$, and
$|z| < (\delta N)^{1/N}$ :
\begin{align}
\| p_K\|_\infty&\leq \| K\| N \leq N\ , \notag\\
\| q_{\del K} \|_\infty &\leq 2\del^2\norm K\norm^2N^2\leq \del N ,\\
\| \Epm^\del \|_\infty&\leq \del N +\del N +\del N\leq
3\del N \notag \ .
\end{align}
Moreover, using also the second line of (\ref{Ebounds}),
\begin{align}
| p_K(0)|& =  | K_{N,1}|, \notag\\
| q_{\del K}(0)| &\leq 2\del^2\norm K\norm^2\leq 2
\del^{2} , \\
| \Epm^\del (0) | &\geq \del( | K_{N,1}| - 2 \delta) \ .
\notag
\end{align}
\end{proof}

To estimate the number of zeros of the holomorphic
function $\Epm^\del(z)$ we will need the next proposition.

\subsection{Counting the zeros and proof of Theorem \ref{th1}}
\begin{proposition}[The Poisson-Jensen formula]\label{jensen}
Let $f$ be a holomorphic function that does not vanish
anywhere on the boundary of $D(0,R)$, where $0<R<\infty$.
Let $M$ be the number of zeros of $f$ in
$D(0,Re^{-\sigma})$ for some positive constant $\sigma$.
Then
\begin{equation}
M \leq \frac{1}{\sigma}\left( - \ln
\frac{|f(0)|}{\|f\|_{L^{\infty}(D(0,R))} } \right) .
\end{equation}
\end{proposition}
This is a direct consequence of formula (1.2'), p.163 in
\cite{neva}: if $f$ is a holomorphic function in $D(0,R)$
with zeros $a_{\mu}$ there, then
\begin{equation}
\ln|f(0)|= \frac{1}{2 \pi} \int_{0}^{2\pi}\ln
|f(Re^{i\theta})| d\theta - \sum_{|a_{\mu}|<R}\ln
\frac{R}{|a_{\mu}|} \leq \ln \|f\|_{\infty}  -
\sum_{|a_{\mu}|<R}\ln \frac{R}{|a_{\mu}|} .
\end{equation}
Hence
\begin{equation}
-\ln |f(0)| + \ln \|f\|_{\infty}\geq \sum_{|a_{\mu}|<R}\ln
\frac{R}{|a_{\mu}|} \geq \sum_{|a_{\mu}|<Re^{-\sigma}}\ln
\frac{R}{Re^{-\sigma}} = M\sigma ,
\end{equation}
which is our proposition.

\begin{proof}[Proof of Theorem \ref{th1}]
With probability at least $1-p_{1}(N)$, we know that
$\|K\|<1$, which we will assume from now on. Then Lemma
\ref{Epm} holds. Setting $f(z)=\Epm^{\delta}(z)$ and
$R=(\delta N)^{1/N}$, this implies that
$\|f\|_{L^{\infty}(D(0,R))}\leq 3\delta N$. We also assume
that $f$ does not vanish on $|z|=R$, since otherwise we
may diminish $R$ slightly and obtain our result as a limit
over increasing radii fulfilling this assumption.

If $\alpha > \delta$, with probability at least $1-
p_{2}(5 \alpha;N)$, we have $|K_{N,1}| - 2 \delta \geq 3
\alpha$, hence again by Lemma \ref{Epm}, $|f(0)| \geq 3
\delta \alpha$.

Thus  we obtain that with probability at least $1-p_{1}(N)
- p_{2}(5 \alpha;N)$,
\begin{equation}
- \ln \frac{|f(0)|}{\|f\|_{L^{\infty}(D(0,R))} } \leq -
\ln \frac{3\delta \alpha}{3\delta N } \leq  \ln N  - \ln
\alpha \ .
\end{equation}
Lemma \ref{Epm} and Proposition \ref{jensen} imply that
with the same probability
\[
\# \big( \Spec(J+\delta K) \cap D(0,(\delta
N)^{1/N} e^{-\sigma }) \big) \leq
\frac{1}{\sigma}(\ln N - \ln \alpha) \ ,
\]
which finishes the proof of Theorem \ref{th1}.
\end{proof}

\subsection{Proof of Theorem~\ref{cor1}}
We end this section by showing that a complex gaussian
random perturbation fulfills the assumptions of Theorem
\ref{th1}.
\begin{proof}
Let $K$ be as in Theorem~\ref{cor1}. Let us first show
that (\ref{normest}) is fulfilled with $p_{1}(N)=
1/N^{2}$.

If $\tilde{K}$ is a random matrix with independent complex
gaussian normal distributed entries, and $a>0$ then
\begin{equation} \begin{split}
P[\|\tilde{K}\|> a]
& \leq P[ \sum_{j,k=1}^N|\tilde{K}_{jk}|^{2} > a^{2}]  \\
&\leq E[a^{-2}\sum_{j,k=1}^N |\tilde{K}_{jk}|^{2}] =
N^{2}/a^{2} \ .
\end{split}\label{eqKsize}
\end{equation}
Hence
\begin{equation} \label{est1}
P[\| K\|< 1] = P[\| \tilde{K}\|< N^{2}]\geq 1- N^{-2}.
\end{equation}
Next, we have to estimate the following probability:
\begin{eqnarray*}
P[|K_{N,1}| \leq s ] &=& P[|\tilde{K}_{N,1}| \leq sN^{2}]\\
& =& 1-\exp\left(
-\frac{(sN^{2})^{2}}{2}\right)\\
&\leq & (sN^{2})^{2} = p_{2}(s;N).
\end{eqnarray*}
Let us choose $\alpha = N^{-3}$. Then
$p_{2}(5\alpha;N)=(5/N)^{2}$, and Theorem \ref{th1}
implies that with probability at least $1-1/N^{2} -
(5/N)^{2}$, we have
\[
\# \big( \Spec(J+\delta K) \cap D(0,(\delta
N)^{1/N} e^{-\sigma }) \big) \leq \frac{4}{\sigma
} \ln N \ ,
\]
which finishes the proof of Theorem~\ref{cor1}.
\end{proof}

\section{Small rank perturbations}
\label{dav}
\subsection{Description of the Example}

Let $C$ be a $k\times k$ matrix with entries $c_{r,s}$ and
let $\del>0$ be the perturbation parameter. We consider
the spectrum of $A:=J+\del K$ where the $N\times N$ matrix
$K$ has the block form
\[
K:=\left( \begin{array}{cc}%
0&0\\
C&0
\end{array}\right)
\]
and the (zero) top right hand entry is of size
$(N-k)\times (N-k)$.

The asymptotic behaviour of $\Spec(A)$ depends on the
choice of $\del$. From this point onwards we assume that
$\del:=\gam^N$ for some $\gam \in (0,\infty)$. If $0<\gam
<1$ then $\del K$ is a very small perturbation of $J$, but
for $\gam >1$ the reverse holds. Surprisingly the same
analysis applies in both cases; the choice $\gam : =1$ is
not special in any way. Our results may have connections
with the analysis of paraorthogonal polynomials on the
unit circle in \cite{Sim}.

The spectral behaviour is different for random
perturbations of the type considered in section \ref{hag}.
The condition number of the diagonalizing matrix is much
smaller in the random model, and the analysis is harder.
For a fixed $\gam \in (0,1)$ the appearance of the
spectrum is similar, but as $\gam \to 1$ the number of
eigenvalues inside the circle increases rapidly, and for
$\gam >1$ the eigenvalues appear to be randomly
distributed.

\begin{example}\label{first} Computations of the eigenvalues
in our model tend to be numerically unstable because of
the high condition numbers involved. We consider the
example in which $N=80$, $\gam =0.6$, $k=3$ and
\[
C:=\left( \begin{array}{ccc} %
8&0&0\\
2&5&0\\
1&-2&3
\end{array}\right)\, .
\]
The eigenvalues of $A$ inside the circle are close to $\pm
i/4$, while the radius of the circle is close to $\gam$.
The condition number of the diagonalizing matrix is
$4.5\times 10^{17}$, and increases for smaller $\gam $ or
larger $N$.
\end{example}

\begin{figure}
\begin{center}
\includegraphics[width=10cm]{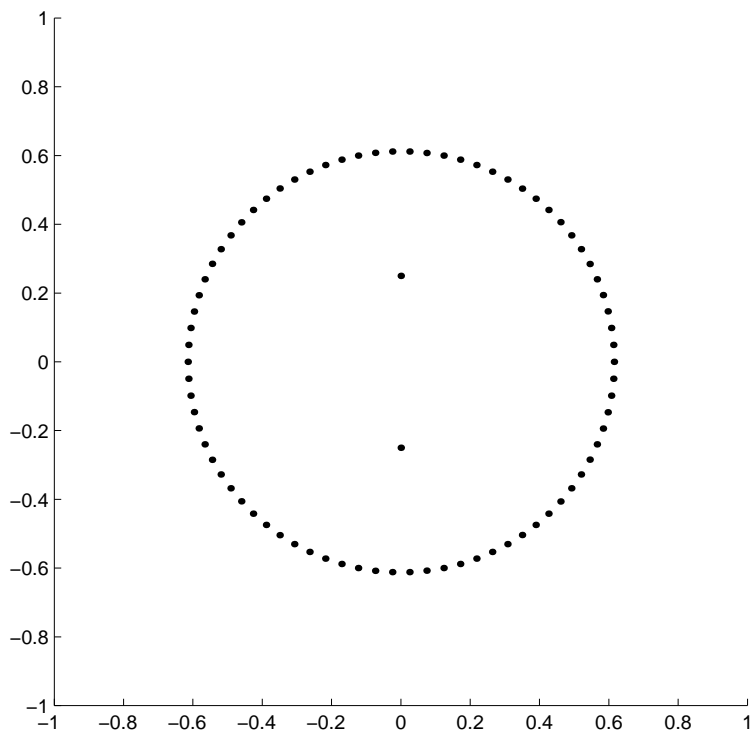}
\end{center}
\caption{Eigenvalues of the Matrix $A$ of
Example~\ref{first}} \label{Figure1}
\end{figure}

\begin{theorem}\label{setup} If $\del:=\gam^N$ where
$\gam \in (0,\infty)$ then
\[
\Spec(A)=\{ \gam z:z^N=f(z)\}
\]
provided $N$ is large enough, where
\[
f(z):=\sum_{r=0}^{2k-2} b_r z^r
\]
and
\[
b_r:=\gam^r\sum_{k-i+j=r+1} c_{i,j}.
\]
\end{theorem}

\begin{proof} The spectrum of $A$ is the set of solutions of
$g_N(\lam)=0$, where
\begin{equation}
g_N(\lam):=\det(\lam I-J-\del K).\label{det}
\end{equation}
Let $g_r(\lam)$ be the determinant of the $r\times r$
matrix obtained from $\lam I-J+\del K$ by deleting its top
$N-r$ rows and the leftmost $N-r$ columns. By expanding
the determinant (\ref{det}) down the leftmost column and
assuming that $N>2(k+1)$ we obtain
\[
g_N(\lam)=\lam g_{N-1}(\lam)- \del c_{k,1}-\del
c_{k-1,1}\lam-\del c_{k-2,1}\lam^2-...-\del
c_{1,1}\lam^{k-1}.
\]
The formula for $g_N(\lam)$ follows inductively, and the
proof is completed by making the change of variables
$\lam:=\gam z$.
\end{proof}

\begin{proof}[Alternative proof using the Grushin problem]
Inserting the special form of $K$ in the series expansion
for $E_{-+}^{\delta}$ and using $N>2(k+1)$, we see that
the series only contains terms up to first order in
$\delta$, so no condition on the smallness of $\gamma$ is
needed for convergence. Putting $\del=\gam^N$ and
$\lam=\gam z$ yields
\begin{equation}
E_{-+}^{\delta}(\lam) = \lam^{N} -\delta E_{-}KE_{+} =
\gamma^{N} ( z^{N} - p_{K}(\gamma z))
\end{equation}
where $p_{K}$ was defined in (\ref{p}). We finally observe
that $p_{K}(\gam z)=f(z)$ for all $z$. Although the
Grushin problem is equivalent to a direct analysis of the
determinant in this particular case, it also permits
estimates in cases where the determinant is quite hard to
analyze directly.
\end{proof}

\subsection{The Equation $z^N=f(z)$}

Let $U$ be a region in the complex plane that contains
$D(0,1+\del)$ for some $\del >0$. Let $f$ be a bounded
analytic function defined on $U$. We assume that $f(z)=0$
has $h$ distinct solutions $z_i$ satisfying $|z_i|<1$,
each with multiplicity $m_i$. We put $n:=\sum_{i=1}^h
m_i$. By reducing $\del >0$ we may assume that
$|z_i|<1-\del$ for all $i$. We will determine the
distribution of the solutions of $z^N=f(z)$ asymptotically
as $N\to\infty$.

\begin{theorem}\label{main1} For every $ \eps \in (0,\del)$ there exists
$N_\eps$ such that if $\, N\geq N_\eps$ then $z^N=f(z)$
has $m_i$ solutions in the $\eps$-neighbourhood of $z_i$
for each $i\in\{1,...,h\}$, no other solutions in
$D(0,1-\eps)$, no solutions in $U\backslash D(0,1+\eps)$
and $N-n$ solutions in $\{ z:1-\eps<|z| <1+\eps\}$.
\end{theorem}

\begin{proof} If $N$ is large enough then $(1+\eps)^N>\max\{|f(z)|:z\in
U\}$, so the equation has no solutions in $U\backslash
D(0,1+\eps)$. By applying Rouche's theorem to $z^N-f(z)$
regarded as a small perturbation of $z^N$, we see that for
all large enough $N$ the equation has $N$ solutions inside
$D(0,1+\eps)$. A similar argument but regarding $f(z)-z^N$
as a small perturbation of $f(z)$, implies that the
equation has $n$ solutions inside $D(0,1-\eps)$, provided
$N$ is large enough, and that these converge to the zeros
of $f(z)$ as $N\to\infty$. The remaining $N-n$ solutions
must lie in the stated annulus.

In order to determine the asymptotic behaviour of the
$N-n$ solutions in the annulus as $N\to\infty$, we assume
for simplicity that $f(z)\not= 0$ whenever $|z|=1$. We
then put
\[
f(\rme^{is}):=\rho(s)\rme^{i\phi(s)}
\]
where $\rho(s)$ is positive and periodic on $[0,2\pi]$
while $\phi(2\pi)=\phi(0)+2\pi n$. Both $\rho(s)$ and
$\phi(s)$ are real analytic functions of $s$. It is easy
to see that for all large enough $N$ the equation
\[
\phi(s)=Ns\hspace{3mm} \mod(2\pi)
\]
has $N-n$ solutions in $[0,2\pi)$. If these are labelled
in increasing order then $s_{r+1}-s_r=2\pi /N+O(1/N^2)$.
We will show that for large enough $N$ the solutions of
$z^N=f(z)$ are very close to the points
$a_r:=\rho(s_r)^{1/N}\rme^{is_r}$.
\end{proof}

\begin{theorem}\label{main2} Given $\alp\in (1,2)$ there
exists a constant $b$ such that for all large enough $N$
and every $r\in \{0,...,N-n-1\}$ the equation $z^N=f(z)$
has a solution $z_r$ satisfying
\[
|z_r-a_r|\leq bN^{-\alp}
\]
To leading order the $N-n$ solutions of $\, z^N=f(z)$ that
are close to the unit circle are uniformly distributed
around it.
\end{theorem}

\begin{proof} An elementary calculation shows that finding
the solution of $z^N=f(z)$ closest to $\rme^{ i s_r}$ is
equivalent to finding the solution of $z^N=f_r(z)$ closest
to $1$, where
\[
f_r(z):=\rme^{- i \phi(s_r)}f(\rme^{ i s_r}z).
\]
We have
\[
f_r(\rme^{is})=\rho_r(s)\rme^{i\phi_r(s)}
\]
where
\begin{eqnarray*}
\phi_r(s)&:=& \phi(s_r+s)-\phi(s_r),\\
\rho_r(s)&:=&\rho(s_r+s).
\end{eqnarray*}
Moreover $\phi_r(0)=0$ and the equation $\phi_r(s)=Ns$ is
equivalent to $\phi(s_r+s)=N(s_r+s)$. From this point
onwards we drop the subscript $r$, assume that
$\phi(0)=0$, and leave the reader to verify that the
bounds obtained are uniform with respect to $r$.

We define the sequence $u_m:=r_m\rme^{i\theta_m}$ for
$m=1,2,...$ by $u_1:=1$ and
$u_{m+1}:=\left\{f(u_m)\right\}^{1/N}$, where we always
take the $N$th root with the smallest argument. Note that
$\theta_1=\theta_2=0$, $r_1=1$ and $r_2=|f(1)|^{1/N}$. In
the following arguments $c_j$ denote positive constants
that do not depend on $N$ or $m$ provided $N$ is large
enough.

We prove that if
\[
S_N:=\{ r\rme^{i\theta}:|r-r_2|\leq N^{-\alp} \mbox{ and }
|\theta|\leq N^{-\alp}\}
\]
then for all large enough $N$, $u\in S_N$ implies
$v:=\{f(u)\}^{1/N}\in S_N$. Put $u:=r\rme^{i\theta}$ and
$v:=s\rme^{i\phi}$. If $u\in S_N$ then
\begin{eqnarray*}
|s^N-r_2^N|&=& |\, |f(u)|-|f(1)|\, | \\
&\leq& |f(u)-f(1)|, \\
&\leq& c_1 |u-1|\\
&\leq& c_1(|u-r_2|+|r_2-1|)\\
&\leq& c_2/N.
\end{eqnarray*}
Therefore
\[
s^N\geq r_2^N-c_2/N=|f(1)|-c_2/N\geq c_3>0
\]
for some $c_3>0$. This also implies that $r_2^N\geq c_3$.
Hence
\[
\sig:=\sum_{i+j=N-1} s^{i}r_2^j\geq Nc_3^{(N-1)/N}\geq
Nc_4.
\]
Combining the above estimates yields
\[
|s-r_2|\leq c_2/N\sig\leq c_5/N^2\leq N^{-\alp}
\]
for all large enough $N$.

We next observe that
\begin{eqnarray*}
Ns^N|\phi|&\leq&\,  c_6|s^N\sin(N\phi)|\,  =\,
c_6|\Im(v^N-r_2^N)|\,  \leq\,
c_6|v^N-r_2^N| \\
&=&c_6|f(u)-f(1)|\,  \leq\,   c_7|u-1|\,  \leq\,   c_8/N.
\end{eqnarray*}
Therefore
\[
|\phi|\leq c_8/c_3N^2\leq N^{-\alp}
\]
for all large enough $N$.

Having established that $S_N$ is invariant under the map
$u\to\{f(u)\}^{1/N}$ provided $N$ is large enough, we now
apply a contraction mapping argument within $S_N$. Let
$z_j\in S_N$ and put
$w_j:=s_j\rme^{i\phi_j}:=\{f(z_j)\}^{1/N}$ for $j=1,2$.
Then
\[
|w_1^N-w_2^N|=|f(z_1)-f(z_2)|\leq c_9 |z_1-z_2|.
\]
Moreover
\begin{eqnarray*}
\rule[-2.1ex]{0.05em}{5ex}\sum_{i+j=N-1}w_1^iw_2^j\,\rule[-2.1ex]{0.05em}{5ex}&\geq
&
\sum_{i+j=N-1}\Re(w_1^iw_2^j)\\
&=& \sum_{i+j=N-1}s_1^is_2^j\cos(i\phi_1+j\phi_2)\\
&\geq & Nc_{10}
\end{eqnarray*}
where $c_{10}>0$. Therefore
\[
|w_1-w_2|\leq c_9|z_1-z_2|/c_{10}N\leq |z_1-z_2|/2
\]
provided $N$ is large enough. Since $u_2\in S_N$, the
contraction mapping principle now implies that the
sequence $u_m$ converges as $m\to\infty$ to a solution
$u\in S_N$ of $u^N=f(u)$, again provided $N$ is large
enough.
\end{proof}

\Note Although we have proved that the eigenvalues of $A$
all lie on or inside the unit circle asymptotically, this
does not imply that $|\det(A)|\leq 1$ asymptotically.
Indeed $\det(A)=(-1)^{N-1}f(0)$ may be of any magnitude.
If $|f(0)|>1$ then the bound
\[
|f(0)|=\left|\frac{1}{2\pi}\int_{-\pi}^\pi
f(\rme^{is})\,\rmd s\right|\leq
\frac{1}{2\pi}\int_{-\pi}^\pi \rho(s)\,\rmd s
\]
implies that $\rho(s)>1$ on average, so the eigenvalues
close to the unit circle are actually slightly outside it,
again on average.

\begin{example}\label{ex8} If there exists $z$ such that $|z|=1$ and $f(z)=0$ then Theorem~\ref{main2} needs to be modified. The estimates in the theorem are local, so the conclusions are applicable to all the solutions of $z^N=f(z)$ that lie in \[ \{z:1-\del\leq |z|\leq 1+\del \mbox{ and } \alp\leq \arg(z)\leq \bet\}, \] provided $f$ does not vanish in this set. Figure~\ref{fignum8} shows the solutions of $z^N=100(z-1)$ when $N=40$.

\begin{figure}
\begin{center}
\includegraphics[width=10cm]{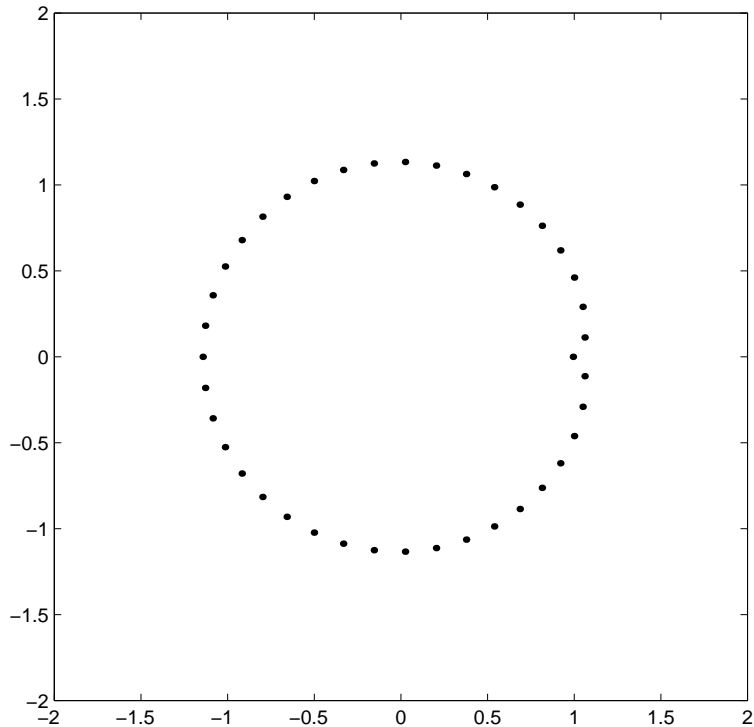}
\end{center}
\caption{Solutions of the polynomial equation of
Example~\ref{ex8}.} \label{fignum8}
\end{figure}
\end{example}

\section{Some Generalizations\label{general}}

\subsection{Other finite rank perturbations}
In this section we allow the perturbation of the Jordan
matrix $J$ to have non-zero entries in all corners of the
matrix. We do not require the perturbation to be small,
since we have already indicated that this possibility can
be accommodated by introducing a scale factor.

Let $B,\, C,\, D,\, E$ be four $k\times k$ matrices and
put
\begin{equation}
A:=J+K\label{genmodel}
\end{equation}
where $N \gg k$ and $K$ has the block form
\[
K:=\left( \begin{array}{ccc}%
B&0&C\\
0&0&0\\
D&0&E
\end{array}\right)\, ,
\]
the central entry being of size $(N-2k)\times (N-2k)$. A
direct calculation shows that
\begin{equation}
\det(\lam I_N-A)=\lam^{N-h}p(\lam)-q(\lam)\label{baseid}
\end{equation}
where $p,\, q$ are polynomials of degree (at most) $2k$
which depend on $B,\, C,\, D,\, E$ but not on $N$, and the
$\lam^{h}$ coefficient of $p$ is $1$.

The following theorem describes the asymptotic
distribution of the eigenvalues of $A$ as $N\to\infty$.
The proof is an obvious adaptation of the proofs of
Theorems~\ref{main1} and \ref{main2}.

\begin{theorem} Let $p,\, q$ be two non-zero polynomials, let
$m$ be a large enough natural number and let $0<\del<1$.
Suppose that $p(z)=0$ and $q(z)=0$ have no solutions
satisfying $1-\del\leq |z|\leq 1+\del$. Then the solutions
of
\[
z^mp(z)=q(z)
\]
satisfying $|z|\leq 1-\del$ converge to the zeros of $q$
in this region as $m\to\infty$. The solutions satisfying
$|z|\geq 1+\del$ converge to the zeros of $p$ in this
region. The solutions satisfying $1-\del\leq |z|\leq
1+\del$ converge to the unit circle and are given
asymptotically by Theorem~\ref{main2}, where
$f(z):=q(z)/p(z)$.
\end{theorem}

\begin{figure}
\begin{center}
\includegraphics[width=10cm]{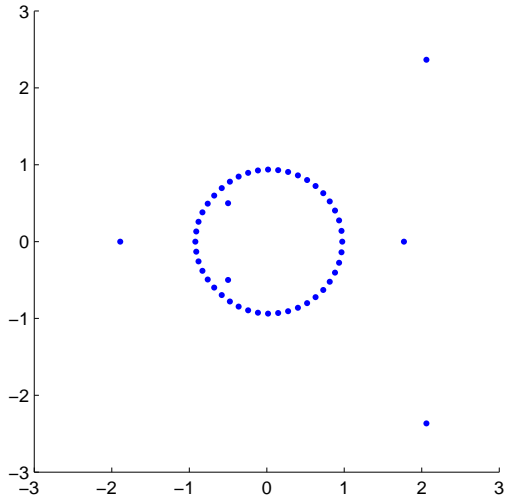}
\end{center}
\caption{Eigenvalues of the matrix $A$ of
Example~\ref{4corners}} \label{Figure2}
\end{figure}

\begin{example}\label{4corners} We consider the above model with
\[
\begin{array}{llll}
\rule[-1.25em]{0em}{4em}B:=&\left(\begin{array}{cc} 2&-1\\ -1&1\end{array}\right), &\hspace{5mm} %
C:=&\left(\begin{array}{cc} 2&-1\\ 1&0\end{array}\right),\\
D:=&\left(\begin{array}{cc} 1&-2\\ 5&3\end{array}\right), & \hspace{5mm}%
E:=&\left(\begin{array}{cc} 1&-3\\ -2&0\end{array}\right).\\
\end{array}
\]
The determinant of $A$ equals $(-1)^N$. Its characteristic
polynomial is $z^{N-4}p(z)-q(z)$, where
\begin{eqnarray*}
p(z)& := & z^4-4z^3+6z^2+15z-33,\\
q(z)& := & -2z^2-2z-1.
\end{eqnarray*}
The zeros of $p$ are all outside the unit circle, at
$2.0605\pm 2.3672i$, $1.7709$ and $-1.8920$. The zeros of
$q$ are both inside the unit circle, at $-1/2\pm i/2$. The
remaining $N-6 $ zeros of the characteristic polynomial
are distributed almost uniformly around the unit circle.
Figure \ref{Figure2} was obtained by putting $N=50$.
\end{example}

\subsection{Matrix Pencils\label{sectpencil}}

A direct calculation shows that the eigenvalues of the
model operator (\ref{genmodel}) are the same as those of
the $(2k+1)\times (2k+1)$ matrix pencil
\begin{equation} \label{an}
A_n(z):= P_n(z)-J_{2k+1}-Q
\end{equation}
where
\[
P_n(z):=\left(\begin{array}{ccc}%
z I_k&0&0\\
0&z^{n-2k}&0\\
0&0&z I_k
\end{array}\right)\, ,%
\hspace{5mm}%
Q:=\left(\begin{array}{ccc}%
B&0&C\\
0&0&0\\
D&0&E
\end{array}\right)\, ,
\]
the central entries of both block matrices being of size
$1\times 1$. By considering other finite rank
perturbations of the Jordan matrix one is led to
investigate the spectral asymptotics as $n\to\infty$ of
the more general $m\times m$ matrix pencil
\[
A_n(z):= B(z)z^n+ C(z)
\]
where $B(z)$  and $C(z)$ are analytic matrix-valued
functions of $z$, defined for all $z$ in an open region
$U$. In this more general problem, $n$ and $A_{n}$ are not
necessarily the same as in (\ref{an}),

The eigenvalues of this pencil are, by definition the
values of $z\in U$ such that $\det(A_n(z))=0$. One sees
that
\[
\det(A_n(z))=\sum_{r=0}^m p_r(z)z^{rn}
\]
where $p_0(z):=\det(C(z))$ and $p_m(z):=\det(B(z))$. If
$p_m$ does not vanish on $U$ then the zeros of the above
expression are the same as those of
\begin{equation}
f_n(z):=z^{mn}+\sum_{r=0}^{m-1} q_r(z)z^{rn}\label{fn}
\end{equation}
where $q_r:=p_r/p_m$ are all analytic functions on $U$.

\begin{theorem} Suppose that the function $f_n$ is defined
by (\ref{fn}) on open region $U$ that contains $B(1)$;
suppose also that $q_r$ are bounded and analytic on $U$
and that $q_0(z)$ does not vanish if $|z|=1$. Then $f_n$
has $mn$ zeros in $U$ for all large enough $n$. They have
the following properties for all small enough $\del >0$.
There are no zeros in $U\cap \{ z:|z|>1+\del\}$ provided
$n$ is large enough. The zeros of $f_n$ satisfying
$|z|\leq 1-\del$ converge as $n\to\infty$ to the zeros of
$q_0$ satisfying the same bound. All the remaining zeros
of $f_n$ lie in $\{z:1-\del\leq |z|\leq 1+\del\}$ and they
converge to the unit circle as $n\to\infty$.
\end{theorem}

The proof is essentially the same as that of
Theorem~\ref{main1}. There is also an analogue of
Theorem~\ref{main2}, but we deal here only with the case
$m=2$. In other words we consider the asymptotic behaviour
of the solutions of an equation of the form
\begin{equation}
p_n(z):= z^{2n}+q_1(z)z^n+q_2(z)=0\label{fundpoly}
\end{equation}
as $n\to\infty$. Examples~\ref{ex2} and \ref{ex3}
illustrate the behaviour that we need to explain.

The following lemma sets up some notation that will be
used in the following theorem.

\begin{lemma}\label{prep} Let $0<\del<1/2$ and let
$q_1$, $q_2$ be two bounded (uniformly in $\delta$)
continuous functions on $A:=\{z:1-\del\leq |z|\leq
1+\del\}$ which are analytic in the interior of this
annulus. Suppose that neither $q_2$ nor the discriminant
$v:=q_1^2-4q_2$ vanish anywhere in $A$. Let $h$, resp.
$k$, be the winding numbers of $v(\rme^{i\theta})$, resp.
$q_2(\rme^{i\theta})$, around the origin, where $\theta\in
[0,2\pi]$. Then the number of solutions of
(\ref{fundpoly}) in $A$ is $2n-k$ for all large enough
$n$.

{\bf Case 1} If $h$ is even then there exist two
non-vanishing analytic functions $f_\pm$ on $A$ such that
$z\in A$ is a solution of (\ref{fundpoly}) if and only if
either $z^n=f_+(z)$ or $z^n=f_-(z)$. Moreover
$q_2(z)=f_+(z)f_-(z)$ for all $z\in A$. If we define the
real-valued analytic functions $\rho_\pm >0$ and
$\phi_\pm$ on $[0,2\pi]$ by
\[
f_\pm(\rme^{i\theta}):=\rho_\pm
(\theta)\rme^{i\phi_\pm(\theta)}
\]
and put $2\pi k_\pm :=\phi_\pm(2\pi)- \phi_\pm(0)$ then
$k:=k_+ +k_-$.

{\bf Case 2} If $h$ is odd then there exists a
non-vanishing analytic function $f$ on the double covering
$\tilde{A}$ of $A$ such that $z\in A$ is a solution of
(\ref{fundpoly}) if and only if $z^n=f(z_1)$ or
$z^n=f(z_2)$, where $z_1$ and $z_2$ are the two points in
$\tilde{A}$ above $z$. Moreover $q_2(z)=f(z_1)f(z_2)$ for
all $z\in A$. If we define the real-valued analytic
functions $\rho
>0$ and $\phi$ on $[0,4\pi]$ by
\[
f(\rme^{i\theta}):=\rho (\theta)\rme^{i\phi(\theta)}
\]
then $2\pi k =\phi(4\pi)- \phi(0)$.
\end{lemma}

\begin{proof} If $|z| \geq 1+\del$, then using the uniform boundedness of $q_{1}$, $q_{2}$ we see that for $n$ big enough the last two terms in (\ref{fundpoly}) are small compared to $|z|^{2n}\gg 1$, hence $p$ cannot vanish there and its $2n$ zeros are confined to a disc:
\[
\frac{1}{2\pi i} \int_{|z|=1+\del}
\frac{p_n^\pr(z)}{p_n(z)}\, \rmd z = 2n \ .
\]
On the other hand, if $|z|=1-\del$, then for $n$ large
enough the first two terms in (\ref{fundpoly}) are so
small that adding them to $q_{2}$ will just have the
effect of moving its zeros slightly inside of $|z|<1-\del$
(recall that these cannot lie on $|z|=1-\del$ by
assumption). Hence in this case they are equal to the
zeros of $p_{n}$ there asymptotically:
\[
\frac{1}{2\pi i} \int_{|z|=1-\del}
\frac{p_n^\pr(z)}{p_n(z)}\, \rmd z = k ,
\]
for all large enough $n$. If $\gam$ is the difference of
these contours then the number of zeros of $p_n$ in $A$ is
given by
\[
\frac{1}{2\pi i} \int_{\gam} \frac{p_n^\pr(z)}{p_n(z)}\,
\rmd z,
\]
which equals  $2n-k$.

To prove the statements in Cases 1 and 2, one only has to
observe that in the formula
\[
2f_\pm(z)= - q_1(z)\pm \sqrt{v(z)}
\]
the square root has a single-valued branch on $A$ if and
only if $h$ is even. The formulae for $k$ follows directly
from $q_2(z)=f_+(z)f_-(z)$ or $q_2(z)=f(z_1)f(z_2)$.
\end{proof}

\begin{theorem} Under the assumptions of Lemma~\ref{prep} the
2n-k solutions of (\ref{fundpoly}) in $A$ converge to the
unit circle. More precisely there exists a constant $c>0$
such that for all large enough $n$ every zero $z\in A$
satisfies
\[
1-c/n\leq |z|\leq 1+c/n.
\]
Moreover the zeros are asymptotically uniformly
distributed around the circle in the sense that
\[
\lim_{n\to\infty} (2n)^{-1}\, \hash \{ z:\theta< \arg
(z)<\phi\} =\phi - \theta
\]
provided $0\leq \theta,\phi\leq 2\pi$.
\end{theorem}

\begin{proof} Case~1 follows directly from Theorem~\ref{main2}, while
Case~2 involves slight modifications of the proof of that
theorem.
\end{proof}

We conclude with two examples exhibiting the behaviour
described in the two cases.

\begin{example}\label{ex2} Consider the equation
\begin{equation}
p_n(z):=z^{2n}+q_1(z)z^n+q_2(z)=0,\label{ex2eq}
\end{equation}
where
\begin{eqnarray*}
q_1(z) &:=& -z^2-z+9/2, \\
q_2(z) &:=& z^3-z^2/2-4z+2.
\end{eqnarray*}
The auxiliary equation
\begin{equation}
w^2+q_1(z)w+q_2(z)=0\label{ex2eq2}
\end{equation}
with $z$ replaced by $\rme^{i\theta}$ has the two distinct
solutions
\begin{eqnarray*}
f_+(\theta) &:=& \rme^{i\theta}-1/2, \\
f_-(\theta) &:=& \rme^{2i\theta}-4,
\end{eqnarray*}
for all $\theta\in [0,2\pi]$. We deduce that  $1/2\leq
|f_+(\theta)|\leq 3/2$ and $3\leq |f_-(\theta)|\leq 5$ for
all $\theta\in [0,2\pi]$. The winding numbers of these
curves around the origin are $m_+:=1$ and $m_-:=0$. The
solutions of $q_2(z)=0$ are $\pm 2$ and $1/2$. We should
therefore anticipate that for large $n$ the equation
(\ref{ex2eq}) has one solution near $z=1/2$ and two
distinct rings of solutions, both close to the unit
circle. One of these rings has $n-1$ points on it while
the other has $n$ points.

This example is particularly simple because one can
factorize (\ref{ex2eq}) in closed form, the left hand side
being the product of $z^n-z+1/2$ and $z^n-z^2+4$. Because
the discriminant of the quadratic equation (\ref{ex2eq2})
is a perfect square its roots come in pairs, so there must
be an even number inside the unit circle. If one replaces
the coefficient $9/2$ of $q_1$ by $9$ one obtains a more
typical example in which the roots form two distinct
rings. The zeros of this modified polynomial
$\tilde{p}_{n}$ are shown in Figure~3, for $n:=40$.
\end{example}

\begin{center}
\scalebox{0.6}{
\includegraphics{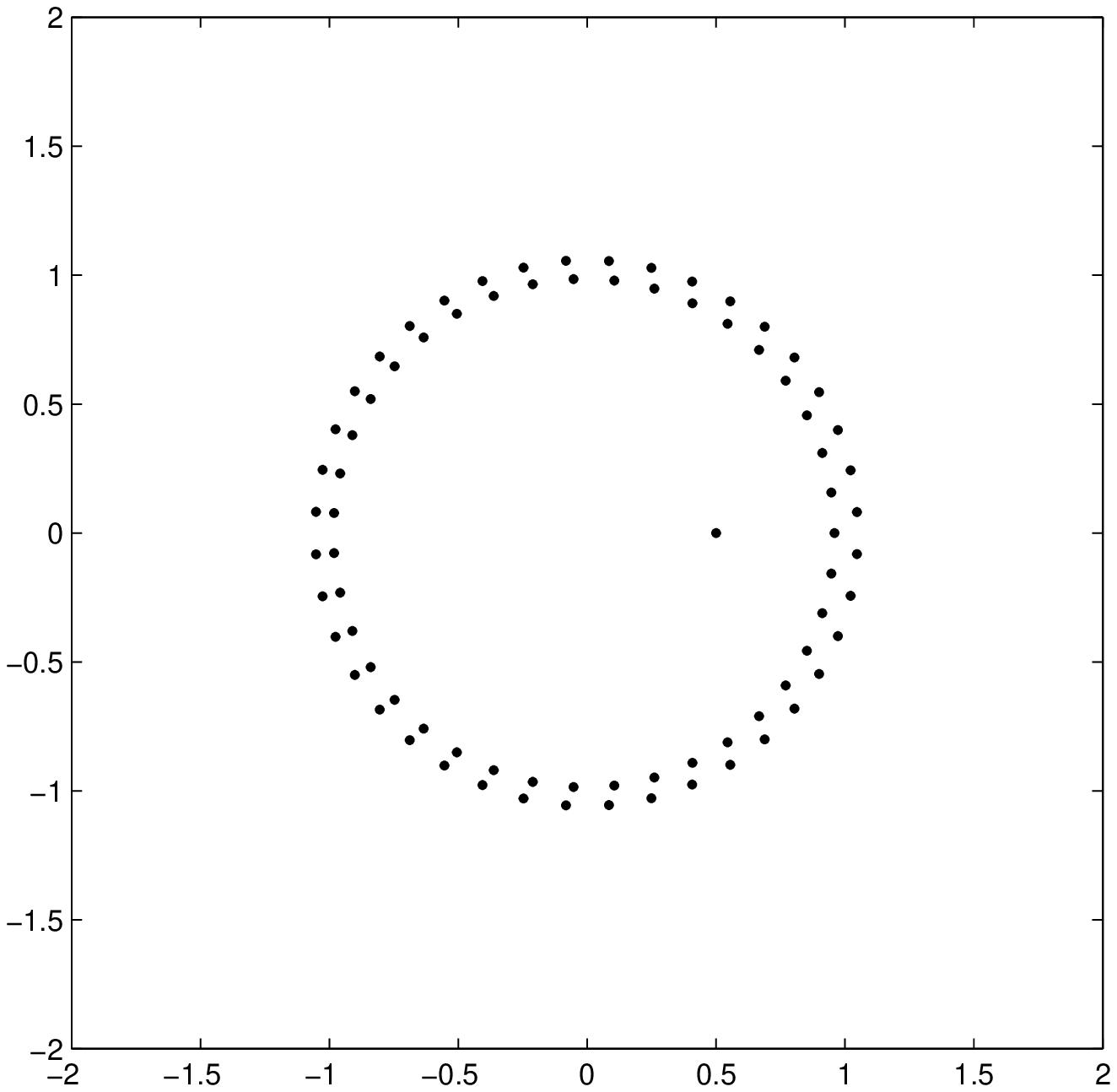}
}\\  %
Figure~3. Eigenvalues of the polynomial $\tilde{p}_{40}$
of Example~\ref{ex2}.
\end{center}

\begin{example}\label{ex3} Consider the equation
\begin{equation}
p_n(z):=z^{2n}-4z^n-8z+3=0.\label{ex3eq}
\end{equation}
The solutions of the auxiliary equation
\begin{equation}
w^2-4w-8\rme^{i\theta}+3=0\label{ex3eq2}
\end{equation}
are
\[
w:=2\pm \sqrt{8\rme^{i\theta}+1}
\]
and combine into a single closed curve winding twice
around the origin and crossing itself on the negative real
axis. Figure~4 shows the set of zeros of the polynomial
equation (\ref{ex3eq}) for $n:=20$.
\end{example}

\begin{center}
\scalebox{0.7}{
\includegraphics{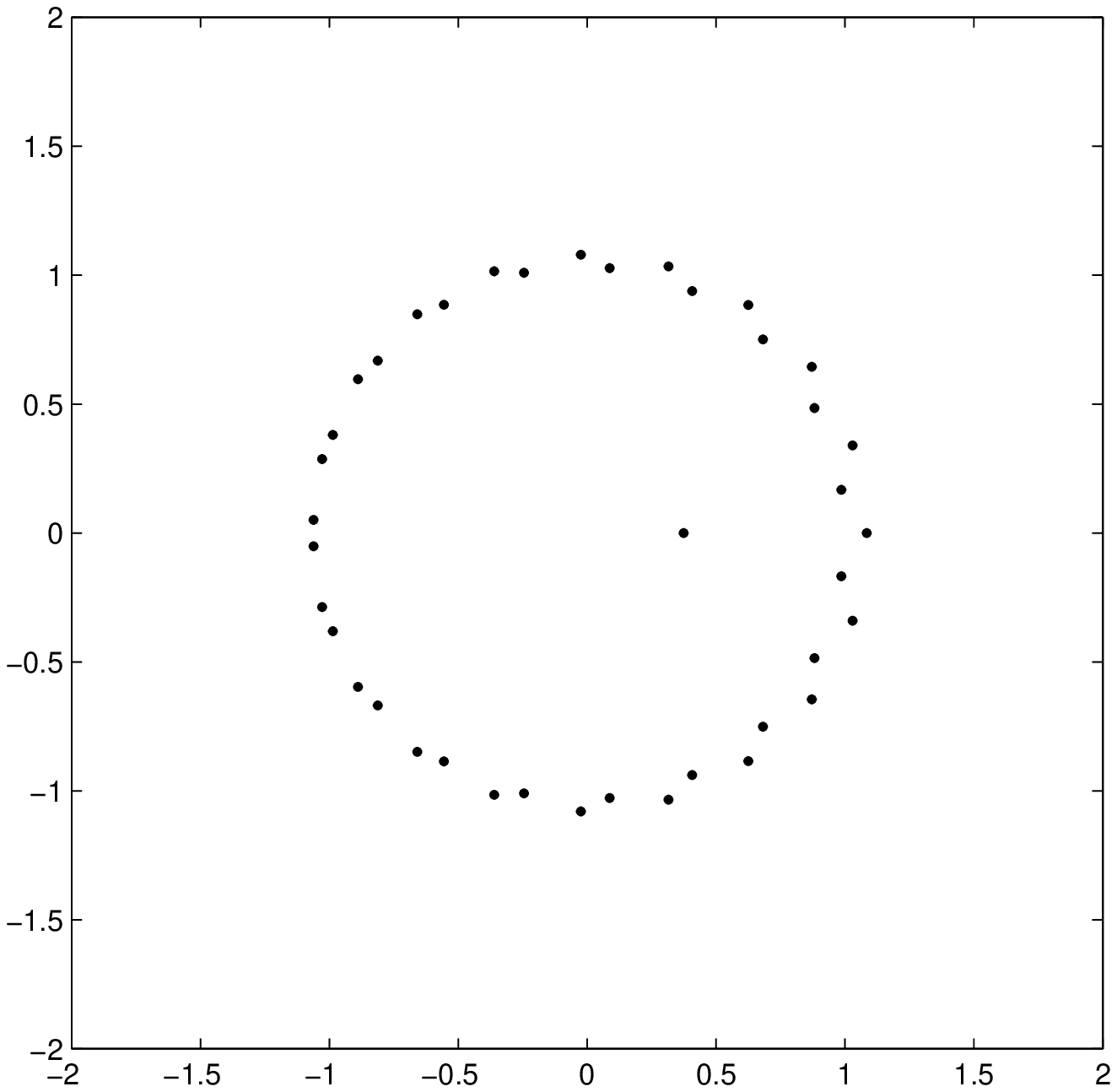}
}\\  %
Figure~4. Eigenvalues of the polynomial $p_{20}$ of
Example~\ref{ex3}.
\end{center}

\medskip

{\bf Acknowledgements} The first author acknowledges
support under EPSRC grant number GR/R81756. The second
author is very grateful to M.~Zworski for pointing out
this problem and its significance to her, and for his
interest in and support of this work. She would also like
to thank S.~Nonnenmacher for his careful reading of a
first draft. She has been supported by a postdoctoral
grant of the Ecole Polytechnique and by a CPAM grant at
the University of California, Berkeley.

\end{document}